\documentclass[11pt, reqno]{amsart}
\usepackage{indentfirst,amssymb,amsmath,amsthm,enumerate, multicol}
\newtheorem*{theoA}{Theorem A}
\newtheorem*{theoB}{Theorem B}
\newtheorem*{theoC}{Theorem C}
\newtheorem*{theoD}{Theorem D}
\newtheorem*{theoE}{Theorem E}
\newtheorem*{theoF}{Theorem F}
\newtheorem*{theoG}{Theorem G}
\newtheorem*{theoH}{Theorem H}
\newtheorem{theo}{Theorem}[section]
\newtheorem{lem}{Lemma}[section]
\newtheorem{cor}{Corollary}[section]

\newcommand{\ol}{\overline}
\newcommand{\be}{\begin{equation}}
\newcommand{\ee}{\end{equation}}
\newcommand{\beas}{\begin{eqnarray*}}
\newcommand{\eeas}{\end{eqnarray*}}
\newcommand{\bea}{\begin{eqnarray}}
\newcommand{\eea}{\end{eqnarray}}
\newcommand{\display}{\displaystyle}

\numberwithin{equation}{section}
\begin{document}
\title[Value sharing by an entire ...]{Value sharing by an entire function with its derivatives}
\author[I.Lahiri \& R. Mukherjee]{Indrajit Lahiri and Rajib Mukherjee }
\begin {abstract}We prove a uniqueness theorem for an entire function, which shares certain values with its higher order derivatives. \end{abstract}
\date{}
\address{Department of Mathematics, University of Kalyani, West Bengal 741235, India.}
\email{ilahiri@hotmail.com}
\address{Department of Mathematics, Krishnanath College, Baharampur, West Bengal 742101, India.}
\email{rajib\_raju786@yahoo.com}
\maketitle
\vspace{-5cm}
\begin{flushleft} Acta Math. Vietnam. (to appear) \end{flushleft}
\vspace{4cm}
\section{Introduction, Definitions and Results}
Let $f$ be a non-constant meromorphic function in the open complex plane $\mathbb{C}$. We denote by $n(r, \infty; f)$ the number of poles of $f$ lying in $\mid z 
\mid < r$, the poles are counted according to their multiplicities. The quantity \[N(r, \infty ; f) = \int\limits_{0}^{r}\frac{n(t, \infty ; f) - n(0, \infty ; 
f)}{t}dt + n(0, \infty ; f) \log r\] is called the integrated counting function or simply the counting function of poles of $f$. 

Also $\display m(r, \infty ; f) = \frac{1}{2\pi}\int\limits_{0}^{2\pi}\log^{+} \mid f(re^{i\theta})\mid d\theta$ is called the proximity function of poles of $f$
, where $\log^{+} x = \log x$ if $x \geq 1$ and $\log^{+} x = 0$ if $0 \leq x < 1$.

The sum $T(r, f) = m(r, \infty ; f) + N(r, \infty ; f)$ is called the Nevanlinna characteristic function of $f$. We denote by $S(r, f)$ any quantity satisfying $S(r, 
f) = o\{T(r, f)\}$ as $r \to \infty$ except possibly a set of finite linear measure.
 
For $a \in \mathbb{C}$, we put $\display N(r, a; f) = N\left(r, \infty ; \frac{1}{f - a}\right)$ and $\display m(r, a; f) = m\left(r, \infty ; \frac{1}{f - 
a}\right)$. 

Let us denote by $\ol n(r, a; f)$ the number of distinct $a$-points of $f$ lying in $\mid z \mid < r$, where $a \in \mathbb{C} \cup \{\infty\}$.
The quantity \[\ol N(r, a ; f) = \int\limits_{0}^{r}\frac{\ol n(t, a ; f) - \ol n(0, a ; f)}{t}dt + \ol n(0, a ; f) \log r\] denotes the 
reduced counting function of $a$-points of $f$.

Also by $\ol N_{(2}(r, a; f)$ we denote the reduced counting function of multiple $a$-points of $f$. 

Let $A \subset \mathbb{C}$ and $n_{A}(r, a; f)$ be the number of $a$-points of $f$ lying in $A \cap \{z : \mid z \mid < r\}$, where $a \in \mathbb{C} \cup 
\{\infty \}$ and the $a$-points are counted acording to their multiplicities. We put 
\[\ol N_{A}(r, a ; f) = \int\limits_{0}^{r}\frac{\ol n_{A}(t, a ; f) - \ol n_{A}(0, a ; f)}{t}dt + \ol n_{A}(0, a ; f) \log r .\]

For $a \in \mathbb{C}\cup \{\infty \}$ we denote by $E(a; f)$ the set of $a$-points of $f$ (counted with multiplicities) and by $\ol E(a; f)$ the set of distinct 
$a$-points of $f$.

For standard definitions and results of the value distribution theory the reader may consult \cite{2} and \cite{6}.

In 1977 L. A. Rubel and C. C. Yang \cite{5} first investigated the uniqueness of entire functions sharing certain values with their derivatives. They proved the 
following result.
\begin{theoA}\cite{5} Let $f$ be a non-constant entire function. If $E(a; f) = E(a; f^{(1)})$ and $E(b; f) = E(b; f^{(1)})$ for two distinct finite complex 
numbers $a$ and $b$, then $f \equiv f^{(1)}$. \end{theoA}

In 1979, E. Mues and N. Steinmetz \cite{4} improved Theorem A in the following manner.
\begin{theoB}\cite{4} Let $a$, $b$ be two distinct finite complex numbers and $f$ be a 
non-constant entire function. If $\ol E(a;f) = \ol E(a;f^{(1)})$ and $\ol E(b;f) = \ol 
E(b;f^{(1)})$, then $f \equiv f^{(1)}$. \end{theoB}

In 1986, G. Jank, E. Mues and L. Volkmann \cite{3} dealt with the case of a single shared value by the two derivatives of an entire function. Their result may be 
stated as follows.
\begin{theoC} \cite{3} Let $f$ be a non-constant entire function and $a (\neq 0)$ be a 
finite complex number. If $\ol E(a; f) = \ol E(a; f^{(1)})$ and $\ol E(a; f) \subset \ol E(a; 
f^{(2)})$ then $f \equiv f^{(1)}$. \end{theoC}

In 2002 J. Chang and M. Fang \cite{1} extended Theorem C in the following way.
\begin{theoD} \cite{1} Let $f$ be a non-constant entire function and $a$, $b$ be two 
non-zero finite constants. If $\ol E(a; f)\subset \ol E(a; f^{(1)}) \subset \ol E(b; f^{(2)})$
, then either $\display f = \lambda e^{\frac{bz}{a}} + \frac{ab - a^{2}}{b}$  or $\display f 
= \lambda e^{\frac{bz}{a}} + a$, where $\lambda (\neq 0)$ is a constant. \end{theoD}

In Theorem C it is not possible to replace the second derivative by any higher order derivative. For, let $f(z) = e^{\omega z} + \omega - 1$, where $\omega^{n - 
1} = 1$, $\omega \neq 1$ and $n (\geq 3)$ is an integer. Then $\ol E(\omega ; f) = \ol E(\omega ; f^{(1)}) = \ol E(\omega ; f^{(n)})$ but $f \not \equiv f^{(1)}$.

Considering higher order derivatives, H. Zhong \cite{8} proved the following result.
\begin{theoE} \cite{8} Let $f$ be a non-constant entire function and $a (\neq 0, \infty)$ be 
a complex number. If $E(a; f) = E(a;f^{(1)})$ and $\ol E(a; f) \subset \ol E(a; f^{(n)}) \cap \ol E(a; f^{(n + 1)})$ for $n \geq 1$, then $f 
\equiv f^{(n)}$. \end{theoE}

P. Li and C. C. Yang \cite{3b} also considered the higher order derivatives and proved the following theorem.
\begin{theoF}\cite{3b} Let $f$ be a non-constant entire function, $a$ be a finite nonzero complex number and $n$ be a positive integer. If $E(a; f) = E(a; 
f^{(n)}) = E(a; f^{(n + 1)})$, then $f \equiv f^{(1)}$. \end{theoF}

To state the next result we require the following definition. Let $f$ and $g$ be two non-constant meromorphic functions defined in $\mathbb{C}$. For $a \in 
\mathbb{C}\cup \{\infty \}$ we put $B = \ol E(a; f) \Delta \ol E(a; g)$, where $\Delta$ denotes the symmetric difference of sets. The functions $f$ and $g$ are 
said to share the value $a$ IMN if $N_{B}(r, a; f) = S(r, f)$ and $N_{B}(r, a; g) = S(r, g)$ \{see \cite{8}\}.

In 1997 L. Z. Yang \cite{7} improved a result of H. Zhong \cite{8} and proved the following theorem.
\begin{theoG}\cite{7} Let $f$ be a non-constant entire function and $a (\neq 0, \infty)$ be a 
complex number. If $f$ and $f^{(n)}$ $(n \geq 1)$ share the value $a$ IMN and $\ol E(a;f) 
\subset \ol E(a; f^{(1)})\cap \ol E(a; f^{(n + 1)})$, then $f = \lambda e^{z}$, where 
$\lambda (\neq 0)$ is a constant.\end{theoG}

Recently Theorem G is improved in the following manner.
\begin{theoH}\cite{3a} Let $f$ be a non-constant entire function and $a (\neq 0, \infty)$ be 
a complex value. Suppose that $A = \ol E(a;f)\backslash \ol E(a;f^{(n)})$ and $B = \ol E(a; 
f^{(n)})\backslash \{\ol E(a; f^{(1)}) \cap \ol E(a; f^{(n + 1)})\}$. If $N_{A}(r,a;f) + 
N_{B}(r, a; f^{(n)}) = S(r, f)$, then either $f = \lambda e^{z}$ or $f = 
\lambda e^{z} + a$, where $\lambda ( \neq 0)$ is a constant. \end{theoH}

It seems to be an interesting problem to investigate the situation when an entire function $f$ shares a nonzero finite value with three consecutive derivatives 
$f^{(n)}$, $f^{(n + 1)}$ and $f^{(n + 2)}$, where $n \geq 1$. In the paper we prove the following result in this direction.
\begin{theo}\label{t1} Let $f$ be a non-constant entire function, $n ( \geq 1)$ be an integer and $a$, $b$ be two nonzero finite complex numbers. Further suppose 
that $A = \ol E(a; f) \backslash \ol E(b; f^{(n)})$ and $B = \ol E(b; f^{(n)}) \backslash \{\ol E(a; f^{(n + 1)}) \cap \ol E(a; f^{(n + 2)})\}.$ If $N_{A}(r, 
a; f) + N_{B}(r, b; f^{(n)}) + \ol N_{(2}(r, a; f) = S(r, f)$, then $a = b$ and either $f = \alpha e^{z}$ or $f = a + \alpha e^{z}$, where $\alpha ( \neq 
0)$ is a constant. \end{theo}

Putting $A = B = \emptyset$ we get the following corollary.
\begin{cor}\label{c1} Let $f$ be a non-constant entire function, $n ( \geq 1)$ be an integer and $a$, $b$ be two nonzero finite complex numbers. If $\ol E(a; f) 
\subset \ol E(b; f^{(n)}) \subset \ol E(a; f^{(n + 1)}) \cap \ol E(a; f^{(n + 2)})$ and $\ol N_{(2}(r, a; f) = S(r, f)$, then $a = b$ and either $f = \alpha 
e^{z}$ or $f = a + \alpha e^{z}$, where $\alpha (\neq 0)$ is a constant. \end{cor}
\section{Lemmas}
In this section we state necessary lemmas.
\begin{lem}\{p.39 \cite{6}\}\label{l1} Let $f$ be a non-constant meromorphic function in 
$\mathbb{C}$ and $n$ be a positive integer. Then \[N(r, 0; f^{(n)}) \leq N(r, 0; f) + n\ol N(r, \infty; f) + S(r, f).\]
\end{lem}
\begin{lem} \{p.57 \cite{2}\}\label{l2} Let $f$ be a non-constant meromorphic function in 
$\mathbb{C}$ and $a$, $b$ be finite nonzero complex numbers and $n$ be a positive integer. 
Then \beas T(r, f) &\leq &\ol N(r, \infty; f) + N(r, a; f) + \ol N(r, b; f^{(n)})\\ &+& S(r, f).\eeas
\end{lem}
\begin{lem} \{p.47 \cite{2}\}\label{l3} Let $f$ be a non-constant meromorphic function in $\mathbb{C}$ and $a_{1}, a_{2}, a_{3}$ be distinct meromorphic functions 
satisfying $T(r, a_{\nu}) = S(r, f)$ for $\nu = 1, 2, 3$. Then \beas T(r, f) &\leq &\ol N(r, a_{1}; f) + \ol N(r, a_{2}; f) + \ol N(r, a_{3}; f)\\ &+& S(r, f),
\eeas where $\ol N(r, a_{\nu}; f) = \ol N(r, 0; f - a_{\nu})$ for $\nu = 1, 2, 3$. \end{lem}
\section{Proof of Theorem \ref{t1}}
\begin{proof}We denote by $N_{(2}(r, a; f \mid f^{(n)} = b)$ the counting function (counted with multiplicities) of those multiple $a$-points of $f$ which are $b$
-points of $f^{(n)}$. We first note that \beas N_{(2}(r, a; f) &\leq & N_{A}(r, a; f) \\ && + N_{(2}(r, a; f \mid f^{(n)} = b) \\ &\leq & n \ol N_{(2}(r, a ; f) + S(r, 
f)\\ &=& S(r, f).\eeas 

Let $z_{1} \not\in A \cup B$ be a simple $a$-point of $f$. Then in some neighbourhood of $z_{1}$ we get by Taylor's expansion 
\beas f(z) & = & a + f^{(1)}(z_{1})(z - z_{1}) + \cdots + \frac{b}{n!}(z - z_{1})^{n} \\ & + & \frac{a}{(n + 1)!}(z - z_{1})^{n + 1} + \frac{a}{(n + 2)!}(z - 
z_{1})^{n + 2} \\ & + & \frac{f^{(n + 3)}(z_{1})}{(n + 3)!}(z - z_{1})^{n + 3} + O(z - z_{1})^{n + 4},\eeas
and so \beas f^{(n)}(z) &=& b + a(z - z_{1}) + \frac{a}{2!}(z - z_{1})^{2} \\ &+& \frac{f^{(n + 3)}(z_{1})}{3!}(z - z_{1})^{3} + O(z - z_{1})^{4},\eeas
\beas f^{(n + 1)}(z) &=& a + a(z - z_{1}) + \frac{f^{(n + 3)}(z_{1})}{2!}(z - z_{1})^{2} \\ &+& O(z - z_{1})^{3}\eeas and 
\beas f^{(n + 2)}(z) = a + f^{(n + 3)}(z_{1})(z - z_{1}) + O(z - z_{1})^{2}.\eeas We note that $f^{(1)}(z_{1}) \neq 0$.

We put $\display \phi = \frac{f^{(n + 1)} - f^{(n + 2)}}{f - a}$, $\display \psi = \frac{f^{(n + 1)} - f^{(n + 2)}}{f^{(n)} - b}$ and $\display H = \frac{bf^{(n 
+ 1)} - a f^{(n)}}{f - a}$. Then by the hypothesis we see that $T(r, \phi) + T(r, \psi) + T(r, H) = S(r, f)$. Now from above we get 
\be\label{1} \phi(z) = \frac{a - f^{(n + 3)}(z_{1})}{f^{(1)}(z_{1})} + O(z - z_{1}), \ee 
\be\label{2} \psi(z) = 1 - \frac{f^{(n + 3)}(z_{1})}{a} + O(z - z_{1}), \ee and 
\be\label{3} H(z) = \frac{ab - a^{2}}{f^{(1)}(z_{1})} + O(z - z_{1}).\ee

We now consider the following cases. \vspace{12pt} \\
{\bf Case 1.} Let $f^{(n + 1)} \equiv f^{(n + 2)}$. Then on integration we get $f^{(n + 1)}(z) = \alpha e^{z}$, where $\alpha (\neq 0)$ is a constant. By 
successive integration we obtain \be\label{4} f(z) = \alpha e^{z} + P(z) = f^{(n + 1)}(z) + P(z),\ee where $P$ is a polynomial of degree $p ( \leq n)$. 

First we suppose that $P$ is non-constant. Then  by Lemma \ref{l3} we get \be\label{5} T(r, f) = \ol N(r, a; f) + S(r, f).\ee 

Now from (\ref{4}) we see that every $a$-point of $f$, which does not belong to $A\cup B$, is a zero of $P$. This shows that \beas \ol N(r, a; f) & \leq & N(r, 
0; P) + N_{A}(r, a; f) + \ol N_{B}(r, a; f) \\ & \leq & N_{B}(r, b; f^{(n)}) + S(r, f) \\ & = & S(r, f), \eeas which contradicts (\ref{5}). Therefore $P(z) 
\equiv \beta$, a constant. Then from (\ref{4}) we get \be\label{6} f(z) = \alpha e^{z} + \beta\ee and so \be \label{7} f^{(n)}(z) \equiv f^{(n + 1)}(z) \equiv 
f^{(n + 2)}(z) = \alpha e^{z}. \ee

We see that $\ol N(r, b; f^{(n)}) \neq S(r, f)$ and $\ol N(r, a; f^{(n + 1)}) \neq S(r, f)$. So by the hypothesis $\ol E(b; f^{(n)}) \cap \ol E(a; f^{(n + 1)}) 
\neq \emptyset$. Hence from (\ref{7}) we get $a = b$. 

Let $\beta \neq a$. Since $f$ does not assume the values $\beta$ and $\infty$, we see that $\ol N(r, a; f) = T(r, f) + S(r, f)$. Again we have from (\ref{7}) 
$\ol N(r, b; f^{(n)}) \neq S(r, f)$. Since $N_{A}(r, a; f) + N_{B}(r, b; f^{(n)}) = S(r, f)$, we get $\ol E(a; f) \cap \ol E(a; f^{(n + 1)}) \neq \emptyset$. So 
from (\ref{6}) and (\ref{7}) we get $\beta = 0$. Therefore $f = \alpha e^{z}$. The other possibility is $\beta = a$ and so $f = a + \alpha e^{z}$. \vspace{12pt} \\
{\bf Case 2.} Let $f^{(n + 1)} \not \equiv f^{(n + 2)}$. By the hypothesis we get 
\bea \ol N(r, b; f^{(n)}) & \leq & N(r, 1; \frac{f^{(n + 2)}}{f^{(n + 1)}}) + N_{B}(r, b; f^{(n)}) \nonumber \\ & \leq & T(r, \frac{f^{(n + 2)}}{f^{(n + 1)}}) + 
S(r, f) \nonumber \\ & = & \ol N(r, 0; f^{(n + 1)}) + S(r, f). \label{8} \eea 

By Lemma \ref{l1} we get from (\ref{8}) \be\label{9} \ol N(r, b; f^{(n)}) \leq N(r, 0; f^{(n)}) + S(r, f).\ee

On the other hand, \beas m(r, a; f) & \leq & m(r, 0; f^{(n)}) + S(r, f) \\ & = & T(r, f^{(n)}) - N(r, 0; f^{(n)}) + S(r, f) \\ &\leq & T(r, f) - N(r, 0; f^{(n)}) 
+ S(r, f) \eeas and so \be\label{10} N(r, 0; f^{(n)}) \leq N(r, a; f) + S(r, f).\ee 

Since $N_{B}(r, b; f^{(n)}) = S(r, f)$, we have $N(r, b; f^{(n)}) = \ol N(r, b; f^{(n)}) + S(r, f)$ and so from (\ref{9}) and (\ref{10}) we get, because $N_{A}(r, 
a; f) = S(r, f)$, \be\label{11} N(r, a; f) = N(r, b; f^{(n)}) + S(r, f).\ee

By Lemma \ref{l2} we obtain from (\ref{11}) \be\label{12} T(r, f) \leq 2 N(r, a; f) + S(r, f).\ee 

First we suppose that $a \neq b$. We put $\display L = \phi - \frac{\psi H}{b - a}$. Then $T(r, L) = S(r, f)$. If possible, let $L \equiv 0$. Then we get $f^{(n 
+ 1)} - f^{(n)} = a - b$. Solving the differential equation we get $f(z) = \alpha e^{z} + P(z),$ where $P$ is a polynomial of degree $n$ with leading coefficient 
$\frac{b - a}{n!}$ and $\alpha$ is a constant. By the hypothesis we see that $f$ cannot be a polynomial and so $\alpha \neq 0$. 

Since $P$ is non-constant, by Lemma \ref{l3} we get \be \label{13} \ol N(r, a; f) = T(r, f) + S(r, f).\ee

Since $N_{A}(r, a; f) = S(r, f)$, by (\ref{13}) we get $\ol E(a; f) \cap \ol E(b; f^{(n)}) \neq \emptyset$. If $z_{0} \in \ol E(a; f) \cap \ol E(b; f^{(n)})$, we 
see that $P(z_{0}) = 0$. Therefore, from (\ref{13}) we get \beas T(r, f) & = & \ol N(r, a; f) + S(r, f) \\ & \leq & N_{A}(r, a; f) + N(r, 0; P) + S(r, f) \\ & = 
& S(r, f),\eeas a contradiction. Hence $L \not \equiv 0$.

Let $z_{1}$ be a simple $a$-point of $f$ such that $z_{1} \not\in A \cup B$. Then by (\ref{1}), (\ref{2}) and (\ref{3}) we get $L(z_{1}) = 0$. Therefore
\beas \ol N(r, a; f) &\leq & N_{A}(r, a; f) + \ol N_{B}(r, a; f) + N(r, 0; L) \\ && +  N_{(2}(r, a; f) \\ & \leq & N_{B}(r, b; f^{(n)}) + S(r, f) \\ & = & S(r, 
f),\eeas which contradicts (\ref{13}). Therefore $a = b$. 

Let $H \not\equiv 0$. If $z_{1} \not \in A \cup B$ is a simple $a$-point of $f$, then from (\ref{3}) we get $H(z_{1}) = 0$. Hence 
\beas \ol N(r, a; f) &\leq & N_{A}(r, a; f) + \ol N_{B}(r, a; f) + N(r, 0; H) \\ && +  N_{(2}(r, a; f) \\ & \leq & N_{B}(r, b; f^{(n)}) + S(r, f) \\ & = & S(r, 
f),\eeas and so $N(r, a; f) \leq \ol N(r, a; f) + N_{(2}(r, a; f) = S(r, f)$, which contradicts (\ref{12}). Therefore $H \equiv 0$ and so $f^{(n)} \equiv f^{(n + 
1)}$. This implies $f^{(n + 1)} \equiv f^{(n + 2)}$, which contradicts the basic assumption of Case 2. This proves the theorem. \end{proof}
\begin{center} {\sc Acknowledgement} \end{center}
The authors are thankful to the referee for his/her valuable suggestions towards the improvement of the paper.
\makeatletter \renewcommand {\@biblabel}[1]{#1.}\makeatother


\begin{thebibliography}{99}
\bibitem{1} J. Chang and M. Fang, Uniqueness of entire functions and fixed points, Kodai 
Math. J., 25 (2002), pp. 309-320.
\bibitem{2} W.K.Hayman, Meromorphic Functions, The Clarendon Press, Oxford (1964).
\bibitem{3} G. Jang, E. Mues and L. Volkmann, Meromorphe Functionen, die mit ihrer ersten 
und zweiten Ableitung einen endlichen Wert teilen, Complex Var. Theory Appl., 6 (1986), pp. 51-71.
\bibitem{3a} I. Lahiri and G. K. Ghosh, Entire functions sharing values with their derivatives, Analysis (Munich), 31 (2011), pp. 47 - 59.
\bibitem{3b} P. Li and C. C. Yang, Uniqueness theorems on entire functions and their derivatives, J. Math. Anal. Appl., 253 (2001), pp. 50 - 57. 
\bibitem{4} E. Mues and N. Steinmetz, Meromorphe functionen, die mit ihrer ableitung werte 
teilen, Manuscripta Math., 29 (1979), pp. 195-206.
\bibitem{5} L. A. Rubel and C. C. Yang, Values shared by an entire function and its 
derivative, in ``Complex Analysis, Kentucky, 1976'', Lecture Notes in Math., Vol. 599, 
Springer (1977), pp. 101-103.
\bibitem{6} C. C. Yang and H. X. Yi, Uniqueness Theory of Meromorphic Functions, Science 
Press and Kluwer Academic Publishers (2003).
\bibitem{7} L. Z. Yang, Further results of entire functions that share one value with their 
derivatives, J. Math. Anal. Appl., 212 (1997), pp. 529-536.
\bibitem{8} H. Zhong, Entire functions that share one value with their derivatives, Kodai 
Math. J., 18 (1995), pp. 250-259.
\end{thebibliography}
\end{document}